\documentclass[10pt]{amsart}
\usepackage{mathrsfs}
\usepackage{amsfonts} %%% i.e. use 12pt type
\usepackage{graphicx,latexsym,bm,amsmath,amssymb,verbatim,multicol,lscape}
\theoremstyle{definition}

\theoremstyle{remark}

\numberwithin{equation}{section}
\begin{document}
\title{The least common multiple of sequence of product of linear polynomials}%
\author{Shaofang Hong}
\address{Mathematical College, Sichuan University, Chengdu 610064, P.R. China}
%\curraddr{}
\email{sfhong@scu.edu.cn, s-f.hong@tom.com, hongsf02@yahoo.com }
\author{Guoyou Qian}
\address{Mathematical College, Sichuan University, Chengdu 610064, P.R. China}
\email{qiangy1230@gmail.com}
\author{Qianrong Tan}
\address{School of Computer Science and Technology, Panzhihua University,
Panzhihua 617000, P.R. China}
\email{tqrmei6@126.com}
\thanks{The research was supported partially by National Science Foundation
of China Grant \# 10971145 and by the Ph.D. Programs Foundation of Ministry
of Education of China Grant \#20100181110073}
\date{\today}%
%\dedicatory{}%
%\commby{}%
% ----------------------------------------------------------------
\begin{abstract}
Let $f(x)$ be the product of several linear polynomials with integer
coefficients. In this paper, we obtain the estimate:  $\log {\rm
lcm}(f(1), ..., f(n))\sim An$ as $n\rightarrow\infty $, where $A$
is a constant depending on $f$.
\end{abstract}

\maketitle

\section{\bf Introduction}

The first significant attempt for proving the prime number theorem
was made by Chebyshev in 1848-1852. In fact, Chebyshev \cite{[Ch]}
introduced the following two functions:
$$\vartheta(x):=\sum_{\ p\le x}\log p=\log \prod_{p\le x}p\
 \ {\rm and}\ \ \psi(x):=\sum_{p^k\le x}\log p=\log
 {\rm lcm}_{1\le i\le \lfloor x\rfloor}\{i\},$$
where $p$ denotes a prime number and  $\lfloor x\rfloor$ denotes the
greatest integer no more than $x$. The prime number theorem asserts
that $\vartheta(n)\sim \psi(n)\sim n$. Thus the asymptotic formula
$\log {\rm lcm}_{1\le i\le n}\{i\}\sim n$ is equivalent to the prime number theorem.
From then on, the topic of estimating the least
common multiple of any given sequence of positive integers become
prevalent and important. Hanson \cite{[Ha]} and Nair \cite{[N]} got the upper bound
and lower bound of ${\rm lcm}_{1\le i\le n}\{i\}$ respectively.
Farhi \cite{[F]} investigated the least common multiple of arithmetic
progression while Farhi and Kane \cite{[FK]} and Hong and Yang \cite{[HY]}
studied the least common multiple of consecutive positive integers.
Recently, Hong and Qian \cite{[HQ]} got some results on the least common multiple of
consecutive arithmetic progression terms. In 2002, Bateman, Kalb
and Stenger \cite{[BKS]} proved that for any integers $a$ and $b$
such that $a\ge 1$ and $a+b\ge 1$ and ${\rm gcd}(a, b)=1$, one has
$$
\log{\rm lcm}_{1\le i\le n}\{ai+b\}\sim \frac{an}{\varphi(a)}\sum_{r=1
\atop \gcd(r, a)=1}^a\frac{1}{r}
$$
as $n\rightarrow\infty $, where $\varphi(a)$ denotes
the number of integers relatively prime to $a$ between 1 and $a$.

Let $h$ and $l$ be any two relatively prime positive integers. The
renowned Dirichlet's theorem says that there are infinitely many
prime numbers in the arithmetic progression $\{hm+l\}_{m\in
\mathbb{N}}$. Furthermore, if we define
$$
\vartheta(x;h,l):=\sum_{{\rm prime}\ p\le x \atop p\equiv l\pmod
{h}}\log p,
$$
then the prime number theorem for arithmetic
progressions says that
$$\vartheta(x;h,l)=\frac{x}{\varphi(h)}+o(x),\eqno (1.1)
$$
For an analytic proof, see Davenport \cite{[D]}. Moreover, Selberg
\cite{[S]} gave an elementary proof of this result.

In this paper, we concentrate on the least common multiple of sequence
of product of linear polynomials with integer coefficients.
Throughout this paper, for any polynomial $g(x)=a_nx^n+\ldots+a_0$
with integer coefficients, we define
$$L_n(g):={\rm lcm}(g(1),\ldots,g(n)).$$
If $\gcd(a_0,\ldots,a_n)=d$, then
$$\log L_n(g)=\log
(dL_n(g_1))=\log L_n(g_1)+\log d=\log L_n(g_1)+O(1),
$$
where $g_1(x)=\frac{a_n}{d}x_n+\ldots+\frac{a_0}{d}$ is a primitive
polynomial. Thus it suffices to give the estimate for  primitive
polynomials. As usual, let  $\mathbb{Q}$ and $\mathbb{N}$ denote the
field of rational numbers and the set of nonnegative integers.
Define $\mathbb{N}^*:=\mathbb{N}\backslash \{0\}$. For any prime
number $p$, we let $v_{p}$ be the normalized $p$-adic valuation of
$\mathbb{N}^*$, i.e., $v_p(a)=s$ if $p^{s}\parallel a$. For any two
positive integers $a$ and $b$, let $\langle b\rangle_a$ denote the
least nonnegative integer congruent to $b$ modulo $a$ between 0 and
$a-1$. We can now state the main result of this paper.\\

\noindent{\bf Theorem 1.1.}  {\it Let $\{t_i\}_{i=0}^k$ be an
increasing sequence of integers with $t_0=0$, and let
$\{d_j\}_{j=1}^{t_k}$ be a sequence of positive integers such that
$d_1=\cdots=d_{t_1}>d_{t_1+1}=\cdots=d_{t_2}>\cdots>d_{t_{k-1}+1}
=\cdots=d_{t_k}$. Let
$$f(x):=\prod_{i=1}^k\prod_{j=t_{i-1}+1}^{t_i}(a_jx+b_j)^{d_j},
$$
where $a_j, b_j\in \mathbb{N}^*$ and $\gcd(a_j,b_j)=1$ for each
$1\le j\le t_k$ and $a_{j_1}b_{j_2}\ne a_{j_2}b_{j_1}$ for any two
integers $1\le j_1\ne j_2\le t_k$. Then we have
$$
\log L_n(f)\sim \frac{n}{\varphi(q)}\sum_{r=1\atop
\gcd(r,q)=1}^q\sum_{i=1}^k (d_{t_i}-d_{t_{i+1}})\max_{1\le j\le
t_i}\{\frac{a_j}{\langle b_jr\rangle_{a_j}}\} \eqno(1.2)
$$
as $n\rightarrow\infty $, where $q={\rm lcm}_{1\le j\le t_k}\{a_j\}$
and $d_{t_{k+1}}:=0$.}\\

Note that if some $b_j$ in Theorem 1.1 are negative integers, then
(1.2) is still true. Evidently, if one picks $k=t_1=d_{t_1}=1$,
then Theorem 1.1 reduces to the Bateman-Kalb-Stenger theorem \cite{[BKS]}.
The proof of Theorem 1.1 will be given in the second section.

\section{\bf Proof of Theorem 1.1}

In this section, we prove the main theorem. For convenience, in what
follows we let $g_j(x):=a_jx+b_j$ for all $1\le j\le t_k$. Then
$f(x)=\prod_{i=1}^k\prod_{j=t_{i-1}+1}^{t_i}g_j(x)^{d_j}$.
We can now show Theorem 1.1 as follows.\\

{\it Proof of Theorem 1.1.} Let $q={\rm lcm}_{1\le j\le
t_k}\{a_j\}$, and let $R(q)=\{r\in \mathbb{N}^*\mid 1\le r\le q,
\gcd(r, q)=1\}$ be the set of positive integers relatively prime to
$q$ and not exceeding $q$. In the following, we let $P_n(f)$ be the
set of the prime factors of $L_n(f)$, and let $P_{n,1}(f)$ denote
the set of the elements in $P_n(f)$ which divide either ${\rm
lcm}_{1\le j_1\ne j_2\le t_k}\{(a_{j_1}b_{j_2}-a_{j_2}b_{j_1})\}$ or
$q$. Define $P_{n,2}(f):=P_n(f)\backslash P_{n,1}(f)$ to be the
complementary set of $P_{n,1}(f)$ in $P_n(f)$. Obviously we have
$$
L_n(f)=\Big(\prod_{p\in P_{n,1}(f)}p^{v_p(L_n(f))}\Big)\Big(\prod_{
p\in P_{n,2}(f)}p^{v_p(L_n(f))}\Big),
$$
equivalently,
$$
\log L_n(f)=\sum_{p\in P_{n,1}(f) }v_p(L_n(f))\log p+\sum_{p\in
P_{n,2}(f)}v_p(L_n(f))\log p. \eqno (2.1)
$$

We claim that if $p\in P_{n,2}(f)$ and $p|f(m)$ for some positive
integer $m$, then there is a unique integer $j_0$ with $1\le j_0\le
t_k$ such that $p|g_{j_0}(m)$ and $p\nmid g_j(m)$ for all other
integers $j$ between $1$ and $t_k$. In fact, we suppose that
$p|(a_{j_1}m+b_{j_1})$ and $p|(a_{j_2}m+b_{j_2})$ for some positive
integers $j_1$ and $j_2$ with $1\le j_1\ne j_2\le t_k$. Then we have
$$p|a_{j_1}(a_{j_2}m+b_{j_2})-a_{j_2}(a_{j_1}m+b_{j_1})
=a_{j_1}b_{j_2}-a_{j_2}b_{j_1}.
$$
It follows that
$$p|{\rm lcm}_{1\le j\ne l\le t_k}\{a_{j}b_{l}-a_{l}b_{j}\},$$
which means that $p\in P_{n, 1}(f)$. This is a contradiction. So the
claim is proved. Thus for any $p\in P_{n,2}(f)$, by the claim we
have
\begin{align*}
v_p(L_n(f))&=\max_{1\le m\le n}\{v_p(f(m))\}=\max_{1\le m\le
n}(\sum_{j=1}^{t_k}d_jv_p(g_j(m)))
=\max_{1\le m\le n}\max_{1\le j\le t_k}\{d_jv_p(g_j(m))\} \\
&=\max_{1\le j\le t_k}\max_{1\le m\le
n}\{d_jv_p(g_j(m))\}=\max_{1\le j\le t_k}\{d_jv_p(L_n(g_j))\}. \ \ \
\ \ \ \ \ \ \ \ \ \ \ \ \ \ \ \ \ \ \ \ (2.2)
\end{align*}

If $p\in P_{n,2}(f)$ and $\max_{1\le j\le t_k}\{v_p(L_n(g_j))\}\ge
2$, then we have $p^2|L_n(g_{j_0})$ for some integer $j_0\in
[1,t_k]$, which implies that $p^2|g_{j_0}(m)$ for some positive
integers $m\le n$. Therefore
$$p\le \sqrt{g_{j_0}(m)}\le
\sqrt{g_{j_0}(n)}\le M_n:=\max_{1\le j\le
t_k}\{\sqrt{g_j(n)}\}.\eqno (2.3)
$$

On the other hand, by the definition of $P_{n,1}(f)$, we obtain that
$P_{n,1}(f)$ consists of only finitely many primes, and hence for
all primes $p\in P_{n,1}(f)$ and all sufficiently large $n$, we have
$p\le M_n\ll\sqrt{n}$. Thus for all sufficiently large $n$, we can
rewrite (2.1) as
 $$
 \log L_n(f)=\sum_{p\le M_n}v_p(L_n(f))\log p+
 \sum_{p>M_n\atop p\in P_{n,2}(f)}v_p(L_n(f))\log p.\eqno (2.4)
 $$

It is obvious that if $p\le M_n$, then
$$v_p(L_n(f))\le \frac{\log f(n)}{\log
p}=\sum_{j=1}^{t_k}d_j\frac{\log g_j(n)}{\log
p}\ll\sum_{j=1}^{t_k}d_j\frac{\log n}{\log p}\ll\frac{\log n}{\log
p}.
$$
Thus we have
 \begin{align*}
\sum_{p\le M_n}v_p(L_n(f))\log p\ll \sum_{p\le M_n}\frac{\log
n}{\log p} \log p &\ll \sum_{p\le M_n}\log n \ll\pi(M_n)\log n\\
& \ll \frac{M_n}{\log M_n}\log n \ll \frac{\sqrt{n}}{\log
\sqrt{n}}\log n\ll \sqrt{n}.
\end{align*}
It then follows from (2.4) that
$$
\log L_n(f)= \sum_{p>M_n\atop p\in P_{n,2}(f)}v_p(L_n(f))\log p
+O(\sqrt{n}).\eqno (2.5)
$$

Now let $p\in P_{n,2}(f)$. Then it is easy to see that $p$ is
congruent to $r'$ modulo $q$ for some $r'\in R(q)$ and
$\gcd(r',a_j)=1$ for all $1\le j\le t_k$. For such $r'$, there is
exactly one $r\in R(q)$ such that $rr'\equiv 1\pmod{q}$, and hence
we have  $\langle b_jr\rangle_{a_j}p\equiv \langle
b_jr\rangle_{a_j}r'\equiv b_jrr'\equiv b_j\pmod{a_j}$ for all $1\le
j\le t_k$. We can easily derive that $ \langle b_jr\rangle_{a_j}p$
is the smallest multiple of $p$ which is congruent to $b_j$ modulo
$a_j$ for all $1\le j\le t_k$. It follows that for any $1\le j\le
t_k$ and any $p\in P_{n,2}(f)$ which is congruent to $r'$ modulo
$q$, we have that $p|(a_jm+b_j)$ for some $m\le n$ if and only if
$p\le \frac{a_jn+b_j}{\langle b_jr\rangle_{a_j}}$. Thus, for any
sufficiently large $n$ and for any prime $p\in P_{n,2}(f)$ which is
congruent to $r'\in R(q)$ modulo $q$ satisfying $p>M_n$, we have by
(2.2) and (2.3) that
$$
v_p(L_n(f))=
\max_{1\le j\le t_k}\{d_jv_p(L_n(g_j))\}=\max_{1\le j\le
t_k}\{e_j\},
$$
where
\begin{align*}
e_j:={\left\{
\begin{array}{rl}
d_j, &\text{if} \ M_n<p\le \frac{a_jn+b_j}{\langle
b_jr\rangle_{a_j}},\\
0, &\text{if}\  p> \frac{a_jn+b_j}{\langle b_jr\rangle_{a_j}}.
\end{array}
\right.}
\end{align*}

Since
$d_1=\cdots=d_{t_1}>d_{t_1+1}=\cdots=d_{t_2}>\cdots>d_{t_{k-1}+1}
=\cdots=d_{t_k}$, we deduce that for sufficiently large $n$,
\begin{align*}
v_p(L_n(f))
={\left\{\begin{array}{rl} d_{t_1}, &\text{if} \ M_n<p\le
\max_{1\le j\le t_1}
\{\frac{a_jn+b_j}{\langle b_jr\rangle_{a_j}}\},\\
d_{t_i}, & {\rm if} \ \max_{1\le j\le
t_{i-1}}\{\frac{a_jn+b_j}{\langle b_jr\rangle_{a_j}}\}< p\le
\max_{1\le j\le t_i}\{\frac{a_jn+b_j}{\langle b_jr\rangle_{a_j}}\} \
\mbox{for\ some} \ 2\le i\le k.
\end{array}
\right.}
\end{align*}
Obviously, we have that for sufficiently large $n$, $v_p(L_n(f))=0$
for any prime $p>M_n$ and $p\not\in P_{n,2}(f)$. Thus we obtain by
(2.5) that
\begin{align*}
&\log L_n(f)=\sum_{p>M_n}v_p(L_n(f))\log p
 +O(\sqrt{n})\\
&=\sum_{r'\in R(q)}\sum_{p>M_n,\atop p\equiv
r'\pmod{q}}v_p(L_n(f))\log p+ O(\sqrt{n})\\
&=\sum_{r'\in R(q)}\Big(\sum_{M_n<p\le \max_{1\le j\le t_1}
\{\frac{a_jn+b_j}{\langle b_jr\rangle_{a_j}}\}\atop p\equiv r'\pmod
q}d_{t_1}\log p\\
& \ \ \ \ \ \ \ \ \ \ \ \ \ \ \ \ \ +\sum_{i=2}^k\sum_{\max_{1\le
j\le t_{i-1}}\{\frac{a_jn+b_j}{\langle b_jr\rangle_{a_j}}\}< p\le
\max_{1\le j\le t_i}\{\frac{a_jn+b_j}{\langle
b_jr\rangle_{a_j}}\}\atop p\equiv r'\pmod q}
d_{t_i}\log p\Big)+O(\sqrt{n})\\
&=\sum_{r'\in R(q)}\Big(d_{t_1}\Big(\vartheta\big(\max_{1\le j\le
t_1}\{\frac{a_jn+b_j}{\langle
b_jr\rangle_{a_j}}\};q,r'\big)-\vartheta(M_n;q,r')\Big)+\sum_{i=2}^k
d_{t_i}F_i(n)\Big)+O(\sqrt{n}),
\end{align*}
where
\begin{align*}
F_i(n):=\vartheta\Big(\max_{1\le j\le
t_i}\big\{\frac{a_jn+b_j}{\langle
b_jr\rangle_{a_j}}\big\};q,r'\Big)-\vartheta\Big(\max_{1\le j\le
t_{i-1}}\{\frac{a_jn+b_j}{\langle b_jr\rangle_{a_j}}\};q,r'\Big).
\end{align*}
Now, applying the prime number theorem for arithmetic progressions
(i.e. (1.1)), we obtain that
\begin{align*}
\log L_n(f)=\frac{n}{\varphi(q)} \sum_{ r'\in R(q)}\sum_{i=1}^k
d_{t_i}\Big(\max_{1\le j\le t_i}\{\frac{a_j}{\langle
b_jr\rangle_{a_j}}\}-\max_{1\le j\le t_{i-1}}\{\frac{a_j}{\langle
b_jr\rangle_{a_j}}\}\Big)+o(n),
\end{align*}
where $\max_{1\le j\le t_{0}}\{\frac{a_j}{\langle
b_jr\rangle_{a_j}}\}:=0$ and $rr'\equiv 1\pmod q$.

Since $r$ runs over $R(q)$ as $r'$ does, it follows immediately that
\begin{align*}
\log L_n(f)&=\frac{n}{\varphi(q)} \sum_{ r\in R(q)}\sum_{i=1}^k
d_{t_i}\Big(\max_{1\le j\le t_i}\{\frac{a_j}{\langle
b_jr\rangle_{a_j}}\}-\max_{1\le j\le t_{i-1}}\{\frac{a_j}{\langle
b_jr\rangle_{a_j}}\}\Big)+o(n)\\
&=\frac{n}{\varphi(q)} \sum_{ r\in R(q)}\Big(d_{t_{k}}\max_{1\le j\le
t_k}\{\frac{a_j}{\langle b_jr\rangle_{a_j}}\}+\sum_{i=1}^{k-1}
(d_{t_i}-d_{t_{i+1}})\max_{1\le j\le t_i}\{\frac{a_j}{\langle
b_jr\rangle_{a_j}}\}\Big)+o(n).
\end{align*}
So we get (1.2) and Theorem 1.1 is proved. \hfill$\Box$\\

In particular, we have the following two consequences.\\

\noindent{\bf Corollary 2.1.}  {\it Let $l\ge 1$ be an integer and
$\{s_i\}_{i=1}^l$ be a decreasing sequence of positive integers, and
let $g(x)=\prod_{i=1}^l(a_ix+b_i)^{s_i}$, where
$a_i, b_i\in \mathbb{N}^*$ and $\gcd(a_i, b_i)=1$ for
each $1\le i\le l$ and $a_ib_j\ne a_jb_i$ for any $1\le
i\ne j\le l$. Then we have
$$\log L_n(g)\sim \frac{n}{\varphi(q)}\sum_{r=1\atop
\gcd(r, q)=1}^q\sum_{i=1}^l(s_{i}-s_{i+1})\max_{1\le j\le
i}\{\frac{a_j}{\langle b_jr\rangle_{a_j}}\},
$$
where $q={\rm lcm}_{1\le i\le l}\{a_i\}$ and $s_{l+1}:=0$.}\\

\noindent{\bf Corollary 2.2.}  {\it Let $l, d\ge 1$ be integers and
$g(x)=\prod_{i=1}^l(a_ix+b_i)^d$, where
$a_i, b_i\in \mathbb{N}^*$ and $\gcd(a_i,b_i)=1$ for
each $1\le i\le l$ and $a_ib_j\ne a_jb_i$ for any two integers $1\le
i\ne j\le l$. Then we have
$$
\log L_n(g)\sim \frac{dn}{\varphi(q)} \sum_{r=1\atop
\gcd(r,q)=1}^q\max_{1\le i\le l}\{\frac{a_i}{\langle
b_ir\rangle_{a_i}}\},
$$
where $q={\rm lcm}_{1\le i\le l}\{a_i\}$. }

\end{document}